\documentclass[12pt]{article}%
\usepackage[latin1]{inputenc}
\usepackage{amsfonts,amssymb,latexsym}
\usepackage[Sonny]{fncychap}
\usepackage[T1]{fontenc}
\usepackage{amsbsy}
\usepackage{graphicx}
\usepackage{times}
\usepackage[latin1]{inputenc}
\usepackage[english]{babel}
\usepackage{amsmath,amsthm}
\usepackage{amsmath}
\usepackage{amsfonts}
\usepackage{amssymb}
\usepackage{color}
\usepackage[colorlinks=true, bookmarksnumbered=true, bookmarksopen=true,
bookmarksopenlevel=3, pdfstartview=FitH, linkcolor=blue, pdfmenubar=true,
pdftoolbar=true, bookmarks=true,citecolor=green, urlcolor=blue,
filecolor=magenta,plainpages=false,pdfpagelabels,breaklinks]{hyperref}%
\newtheorem{teor}{Theorem}[section]

\newtheorem{lem}[teor]{Lemma}
\newtheorem{prop}[teor]{Proposition}
\theoremstyle{definition}
\newtheorem{defn}[teor]{Definition}
\theoremstyle{remark}
\newtheorem{rem}[teor]{Remark}

\numberwithin{equation}{section}

\def\fin { \vskip 0pt \hfill \hbox{\vrule height 5pt width 5pt depth 0pt} \vskip 12pt}
\begin{document}

\title{Periodic solutions for a 1D-model with nonlocal velocity via mass transport}
\author{{{Lucas C. F. Ferreira} {\thanks{L. Ferreira was supported by FAPESP and CNPQ,
Brazil. (corresponding author)}}}\\{\small Universidade Estadual de Campinas, Departamento de Matem\'{a}tica,}\\{\small {\ CEP 13083-859, Campinas-SP, Brazil.}}\\{\small \texttt{E-mail:lcff@ime.unicamp.br}}\vspace{0.5cm}\\{Julio C. Valencia-Guevara }\\{\small Universidade Estadual de Campinas, Departamento de Matem\'{a}tica,}\\{\small {\ CEP 13083-859, Campinas-SP, Brazil.}}\\{\small \texttt{E-mail:ra099814@ime.unicamp.br}}}
\date{}
\maketitle

\begin{abstract}
This paper concerns periodic solutions for a 1D-model with nonlocal velocity
given by the periodic Hilbert transform. There is a rich literature showing
that this model presents singular behavior of solutions via numerics and
mathematical approaches. For instance, they can blow up by forming
mass-concentration. We develop a global well-posedness theory for periodic
measure initial data that allows, in particular, to analyze how the model
evolves from those singularities. Our results are based on periodic mass
transport theory and the abstract gradient flow theory in metric spaces
developed by Ambrosio et al. [2]. A viscous version of the model is also
analyzed and inviscid limit properties are obtained.

\

\

\vspace{0.1cm}


\end{abstract}

\

\section{Introduction}

\hspace{0.4cm} We consider the following one-dimensional model
\begin{equation}
u_{t}+(\mathcal{H}(u)u)_{x}=0 \label{eq-1}%
\end{equation}
with the initial condition $u(x,0)=u_{0},$ where $\mathcal{H}$ stands for the
periodic Hilbert transform
\begin{equation}
\mathcal{H}(u)(x)=\frac{1}{2\pi}\text{P.V.}\int_{-\pi}^{\pi}\cot(\frac{x-y}%
{2})u(y)dy, \label{trans1}%
\end{equation}
and the unknown $u:\mathbb{R}\times\lbrack0,\infty)\rightarrow\mathbb{R}$ is
$2\pi$-periodic in the spatial variables. The viscous version of (\ref{eq-1})
is also studied. In the non-periodic case, (\ref{trans1}) should be changed to
the continuous Hilbert transform
\begin{equation}
\mathcal{H}_{c}(u)(x)=\frac{1}{\pi}\text{P.V.}\int_{\mathbb{R}}\frac
{u(y)}{x-y}{dy}. \label{trans2}%
\end{equation}

The model (\ref{eq-1}) is a continuity equation with non-local velocity field
and has a strong analogy with some physical models; for instance, 2D inviscid
quasi-geostrophic equation (see \cite{Chae-Cordoba-Fontelos}) and 2D vortex
sheet problems (see \cite{Baker-Morlet}). It also appears in modeling of
dislocation dynamics in crystals where $u\geq0$ stands for the density of
defects in the material (see \cite{Deslippe},\cite{Head},\cite{Biler-Karch}).
Applying the Hilbert transform over (\ref{eq-1}), and computing the
$x$-derivative, the resulting equation can be used, in a first approximation,
to study the dynamics of the interface between two fluids; one governed by
Stokes equations and other by Euler equations (see \cite[Appendix
A]{Castro-Cordoba}). For other examples of 1D-equations appearing as models
for PDEs defined in higher dimensions, we refer the reader to
\cite{Constantin1},\cite{Cordoba-Fontelos1},\cite{Escudero},\cite{Gregory1}
and their references. In particular, the non-conservative variant of
(\ref{eq-1})
\begin{equation}
u_{t}+\mathcal{H}(u)u_{x}=0 \label{eq-3}%
\end{equation}
and its viscous versions have been studied by several authors via mathematical
fluid mechanics arguments, see \cite{Cordoba-Fontelos1}%
,\cite{Cordoba-Fontelos2},\cite{Dong-H},\cite{Dong-Li} and their references.

Beyond that, (\ref{eq-1}) has a mathematical interest of its own due to its
non-local structure and singular behavior with respect to existence of global
solutions. For instance, in comparison with (\ref{eq-3}), a difficulty of
handling (\ref{eq-1}) is the lack of maximum principle for the $L^{\infty}%
$-norm. In fact, in \cite{Chae-Cordoba-Fontelos}, the authors showed there is
no global periodic solutions of (\ref{eq-1}) in $C^{1}([-\pi,\pi]\times
\lbrack0,\infty))$ for $u_{0}\in$ $C^{1}([-\pi,\pi])$ with $\int_{-\pi}^{\pi
}u_{0}(x)dx=0$ and $u_{0}\not \equiv 0.$ If, instead, one assumes
\begin{equation}
\int_{-\pi}^{\pi}u_{0}(x)dx\geq0\text{ and }\min_{x}u_{0}(x)<0 \label{cond-1}%
\end{equation}
then the $C^{1}$-breakdown still holds true. The authors of
\cite{Castro-Cordoba} considered the non-periodic version of (\ref{eq-1}) and
showed local well-posedness of nonnegative $H^{2}(\mathbb{R})$-solutions.
These develop a finite time singularity provided that there is a $x_{0}%
\in\mathbb{R}$ such that $u_{0}(x_{0})=0.$ For $0<\delta<1$ and strictly
positive $u_{0}\in L^{2}(\mathbb{R})\cap C^{1,\delta}(\mathbb{R})$ vanishing
at infinity, they showed there is a unique global solution $u\in
C([0,\infty);L^{2}\cap C^{1,\delta})$ for (\ref{eq-1}) and a version of it
with viscous term $-\nu\mathcal{H}_{c}(u_{x}).$ Considering the fractional
viscosity $(-u_{xx})^{\alpha/2},$ the paper \cite{Dong-Li-2} extended some
results in \cite{Castro-Cordoba} by showing blow up of smooth solutions for
$0\leq\alpha\leq2$ and smooth positive initial data with sufficiently
localized mass. The authors of \cite{Baker-Morlet} studied periodic solutions
for
\begin{equation}
u_{t}+(\mathcal{H}(u)u)_{x}=\nu u_{xx} \label{eq1-viscous}%
\end{equation}
which is (\ref{eq-1}) with the additional viscous term $\nu u_{xx}.$ There,
the Hilbert transform is taken with an opposite sign but solutions correspond
easily by changing $u$ to $-u$. For $\nu>0,$ they constructed a solution $u\in
C^{\infty}(\mathbb{T\times}[0,T^{\ast}))$ with $L^{\infty}$-norm blowing up at
a certain finite time $T^{\ast}.$ In fact, they showed that $u$ converges in
$\mathcal{D}^{\prime}(\mathbb{R})$ as $t\rightarrow T^{\ast}$ to the periodic
measure
\begin{equation}
\mu_{a}=a+\Sigma_{n\in\mathbb{Z}}\delta_{0}(x-2\pi n), \label{data-1}%
\end{equation}
where $a$ is a constant and $\delta_{0}$ is the Dirac delta distribution. In
other words, by standard periodic identification, $u(\cdot,t)\rightharpoonup
a+\delta_{0}|_{[-\pi,\pi)}$ in $\mathcal{D}^{\prime}(\mathbb{T})$ as
$t\rightarrow T^{\ast}.$ For $\nu=0$ and a positive data $u_{0}\in C^{\infty
}(\mathbb{T)}$, blow up of $L^{\infty}(\mathbb{T})$-norm was proved in
\cite{Baker-Morlet} which also indicates a concentration of mass due to
sign-preservation and mass-conservation for solutions of (\ref{eq-1}). Taking
data in the form
\begin{equation}
u_{0}(x)=a_{0}+a_{1}\cos(x) \label{data-periodic}%
\end{equation}
where $\left\vert a_{1}\right\vert >\nu\geq0$ and $a_{0}\neq0$, solutions with
$L^{2}(\mathbb{T})$-norm blowing up at a finite time were obtained in
\cite[p.157]{Morlet}. These solutions can (or not) be nonnegative according to
the choice of the parameters $a_{1}$ and $a_{0}$.

The above results corroborate the singular feature of (\ref{eq-1}) and, in
particular, show that solutions can exhibit mass concentration. So, it is
natural to wonder about a framework in which solutions could continue after
the blow-up time and how the PDE evolves from singular data like
(\ref{data-1}). Motivated by that, we investigate (\ref{eq-1}) in a setting of
periodic measures and show global well-posedness of the gradient flow
associated to (\ref{eq-1}). In fact, we consider (\ref{eq1-viscous}) with
$\nu=0$ and $\nu>0$ both for $u_{0}$ belonging to the set of periodic
probability measures $\mathcal{P}(\mathbb{S}^{1})$ endowed with the periodic
Wasserstein metric (see \cite{Ambrosio-2},\cite{Gangbo},\cite{Cordero-1}). For
all initial data $u_{0}\in\mathcal{P}(\mathbb{S}^{1})$, solutions converge
towards a stationary state as $t\rightarrow\infty,$ which is the unique
minimum for the associated energy functional. Moreover, we prove that
solutions of (\ref{eq1-viscous}) with $\nu>0$ converge in $\mathcal{P}%
(\mathbb{S}^{1})$ to those of (\ref{eq-1}) when $\nu\rightarrow0^{+}$(inviscid
limit). In view of the mass-conservation property, notice that (by making a
normalization) the constraint $\int_{-\pi}^{\pi}u_{0}dx=1$ is not an essential one.

We also point out that the evolution of (\ref{eq-1}) from initial measures may
be of interest due to its connection with some problems involving 2D vortex
sheet which is in a layer of vorticity distributed as a delta function on a
curve. In \cite{Chae-Cordoba-Fontelos}, an explicit formula for solutions of
(\ref{eq1-viscous}) with $\nu>0$ was obtained by using the Hopf-Cole transform
and complex Burgers equation. To do this, it is necessary to have
$\mathcal{H}(u_{0})$ at least belonging to $L_{loc}^{1}$ what is not verified
for a general $u_{0}\in\mathcal{P}(\mathbb{S}^{1})$; for instance,
$\mathcal{H}(\delta_{0})\notin L_{loc}^{1}(-\delta,\delta),$ for all
$0<\delta<\pi.$

Formally, the PDEs (\ref{eq-1}) and (\ref{eq1-viscous}) can be rewrite as a
continuity equation
\begin{equation}
\partial_{t}u(x,t)=\nabla\cdot(\upsilon(x,t)u(x,t))\text{ } \label{continuity}%
\end{equation}
with velocity field $\upsilon=\nabla\frac{\delta E}{\delta u}$ given by the
gradient of the variational derivative of the corresponding free energy
functional (see Section 2 for details). Equations in this form have the
so-called gradient-flow structure (see \cite{Ambrosio}) and their solutions
can be obtained by means of an interactive variational scheme based on optimal
transport theory and properties of $E$, what goes back to the seminal work
\cite{Jordan-Otto} for the linear Fokker-Planck equation. Roughly speaking,
the basic idea is to construct solutions that follows the direction of
steepest descent of the energy functional in a probability measure space
endowed with a suitable metric. In the non-periodic setting, an appropriate
space is $\mathcal{P}_{2}(\mathbb{R}^{d})$ (probability measures with finite
second moments) endowed with the $2$-Wasserstein distance.

In fact, the above approach has been applied in $\mathcal{P}_{2}%
(\mathbb{R}^{d})$ to several equations (see \cite{Agueh},\cite{Carrillo-2}%
,\cite{Otto}) and an abstract theory has been developed to a general class of
continuity equation in $\mathbb{R}^{d}$ with energy functional
\begin{equation}
E[u]:=\int_{\mathbb{R}^{d}}U(u(x))\,dx+\int_{\mathbb{R}^{d}}%
u(x)\,V(x)\,dx+\frac{1}{2}\iint_{\mathbb{R}^{d}\times\mathbb{R}^{d}%
}W(x-y)\,u(x)\,u(y)\,dx\,dy, \label{generalfunctionals}%
\end{equation}
where the terms $U:\mathbb{R}^{+}\rightarrow\mathbb{R}$, $V:\mathbb{R}%
^{d}\rightarrow\mathbb{R}$ and $W:\mathbb{R}^{d}\rightarrow\mathbb{R}$ are a
density of internal energy, a confinement potential and an interaction
potential, respectively (see e.g. \cite{Ambrosio},\cite{Carrillo1}). The
concept of displacement convexity for $E$ introduced in \cite{McCann} (more
generally, $\lambda$-convexity) plays a core role in the theory, as well as
lower semicontinuity and coercivity properties of $E$. More precisely, in
\cite[Chapter 4]{Ambrosio}, a gradient-flow theory is developed for
(\ref{continuity}) in general metric spaces by assuming these properties for
an abstract functional $E$. An important example is the Wasserstein metric
space $\mathcal{P}_{2}(X)$ where $X$ is a separable Hilbert space. This
general theory was successfully used in \cite[Chapter 11, p. 298-303]%
{Ambrosio} to study some PDEs in $X=\mathbb{R}^{d}$ with functionals having
the concrete form (\ref{generalfunctionals}) and satisfying certain smooth and
growing conditions on $U,V,W.$

In \cite{Car-Fer-Pre}, the authors analyzed (\ref{eq-1}) and
(\ref{eq1-viscous}) in the non-periodic case where the interaction potential
$W(x)=-\frac{1}{\pi}\log\left\vert x\right\vert $ is singular at origin.
There, by employing the results from \cite[Chapter 11]{Ambrosio}, a
gradient-flow solution in $\mathcal{P}_{2}(\mathbb{R}^{d})$ was obtained after
making a self-similar change of variables in the equations and proving the
needed properties for $E$. This change of variables generates a confinement
term $V(x)=\frac{\left\vert x\right\vert ^{2}}{2}$ that helps to control the
lack of boundedness from below of the interaction part of $E$. The asymptotic
behavior of solutions in \cite{Car-Fer-Pre} is described by a self-similar one while here the dynamics is attracted to a unique stationary solution.

Let us also comment about motivation from a general point of view. Periodic
solutions are widely studied in PDE-theory and appear naturally in several
physical phenomena, specially in fluid mechanics. So, it is important that
different approaches get to deal with this kind of solution. In this direction, our results seem to be the first construction of space-periodic
gradient flows in the context of fluid mechanics. In fact, there are a few
works dealing with existence of periodic solutions for PDEs via optimal mass
transport. For instance, we would like to mention the papers \cite{Carrillo-3}%
, \cite{Ambrosio-2} and \cite{Gangbo}. In the former, the authors analyzed the
family of first-order displacement-convex functionals in $\mathcal{P}%
(\mathbb{S}^{1})$
\begin{equation}
E(\rho)=\int_{0}^{1}\left[  \left(  \frac{1}{\rho^{\beta}}\right)
_{x}\right]  ^{2}dx,\text{ for }\beta\in\lbrack1,\frac{3}{2}], \label{func-1}%
\end{equation}
whose associated gradient flows are periodic weak-solutions of a class of
fourth-order degenerate parabolic equations. In \cite{Ambrosio-2}, existence
of Eulerian distribution solutions for semigeostrophic equations was obtained
by using regularity and stability properties found in \cite{DePhili-Fig-1,
DePhili-Fig-2} for Alexandrov solutions of the Monge-Ampere equation and
optimal mass transport in $\mathbb{T}^{2}$. The authors of \cite{Gangbo}
developed a weak KAM theory in $\mathcal{P}(\mathbb{T}^{d})$ with $d\geq1$ and
obtained existence and asymptotic behavior of solutions for the nonlinear
Vlasov system.

Another motivation is that objects like (\ref{data-1}) and
(\ref{data-periodic}) do not belong to $\mathcal{P}_{2}(\mathbb{R})$ and then
they are not covered by the results in \cite{Car-Fer-Pre}. Moreover, it is
worthy to mention that periodic conditions prevent the use of the self-similar
change employed by \cite{Car-Fer-Pre}. So, we need to handle the singular
interaction potential of (\ref{eq-1}) and (\ref{eq1-viscous}) in original
variables in order to obtain the key properties for $E$ and carry out in
$\mathcal{P}(\mathbb{S}^{1})$ the general theory in metric spaces of
\cite[Chapter 4]{Ambrosio}.

In what follows, we comment on some technical difficulties. Unlike when $X$ is
Hilbert, the sphere $\mathbb{S}^{1}$ is not a convex set and then a
displacement interpolation curve like ${((1-t)H_{1}+tH_{2})_{\#}}\mu$ with
$H_{i}:\mathbb{S}^{1}\rightarrow\mathbb{S}^{1}$ could not be well-defined in
$\mathcal{P}(\mathbb{S}^{1})$. In Section 3.3, we work with a concept of
generalized geodesic in $\mathcal{P}(\mathbb{S}^{1})$ as curves of equivalence
classes (see Definition \ref{Def-Geo-1} and Remark \ref{rem-0}). By using
this, we define and show a type of convexity for functionals. In particular,
in despite of the cost $d_{per}^{2}(x,\cdot)$ in (\ref{dist-1}) is not convex,
we show the $2$-convexity of the square of the periodic Wasserstein distance
$\mathbf{d}_{per}^{2}(\mu,\cdot)$ (see Lemma \ref{convexitywasserstein}),
which is essential for the convergence of the steepest descent scheme
(\ref{scheme-1})-(\ref{discrete-sol-1}). This is obtained by employing the
equivalent representations $\mathcal{P}(\mathbb{S}^{1})$ and $\mathcal{P}%
_{2}(\mathbb{R})/\sim$ (see (\ref{dist-1}) and (\ref{rel-equi1})) and the
identity (\ref{relmetric1}). Also, in Lemma \ref{entropyinvariance}, we prove
a certain invariance property for $\int_{\mathbb{R}^{d}}U(u(x))\,dx$ (where
$U(0)=0$) with respect to the equivalence relation (\ref{rel-equi1}). This is
key in the proof of the convexity of the entropy functional $\mathcal{U}%
[u]=\int_{[-\pi,\pi)}u\log u dx$ insofar as it assures the invariance of the
integral (\ref{aux-func-1}) for elements of an equivalent class in
$\mathcal{P}_{2}(\mathbb{R})/\sim$ supported in some interval of the type
$[a,a+2\pi).$ Similarly, in order to obtain convexity of the interaction
functional $\int\int_{[-\pi,\pi)^{2}}W(x-y)d u(x)d u(y)$, we need an
invariance property that is stated in Remark \ref{welldefinition1}.

Furthermore, to our knowledge, there is no stability result in general metric
spaces for solutions obtained via the theory of gradient flows found in
\cite[Chapter 4]{Ambrosio}. The authors of \cite{Car-Fer-Pre} employed a
stability result of \cite[Chapter 11]{Ambrosio} in the space $\mathcal{P}%
_{2}(X)$ where $X$ is Hilbert. Since this is not the case of $\mathbb{S}^{1}$,
it is necessary to obtain a version of such results for the periodic setting
(see Theorem \ref{lim-inviscid}). In fact, the periodic condition allows to
perform a proof more direct than that in \cite[Chapter 11]{Ambrosio} and could
be extended to study stability of gradient flows in $\mathcal{P}%
(\mathbb{S}^{1})$ generated by a family of general functionals $\{E_{\alpha
}\}$ under relatively simpler conditions.

The plan of this paper is as follows. In Section 2, we summarize some facts
about optimal mass transport in $\mathbb{S}^{1}$ and present the gradient-flow
structure of (\ref{eq-1}) and (\ref{eq1-viscous}) in a more detailed way.
Section 3 is devoted to prove key properties of the free energy functional
$E$. Finally, global well-posedness of gradient-flow solutions in
$\mathcal{P}(\mathbb{S}^{1})$ and inviscid limit are proved in Section 4.

\section{Mass transport and Gradient-Flow Structure}

\subsection{\bigskip Mass Transport in $\mathbb{S}^{1}$}

\hspace{0.4cm} In this section, we resume the theory of optimal transport
relevant for our purposes. The circle $\mathbb{S}^{1}$ is considered as the
quotient space $\mathbb{R}/2\pi\mathbb{Z}$ and functions on $\mathbb{S}^{1}$
are considered as $2\pi$-periodic functions in $\mathbb{R}$. The space
$C^{r}(\mathbb{S}^{1})$ stands for the set of $2\pi$-periodic functions of
class $C^{r}$, for $r\geq0.$ In the case $r=0,$ we denote $C^{0}%
(\mathbb{S}^{1})$ by $C(\mathbb{S}^{1}).$ Also, $\mathbb{S}^{1}$ will be
identified with the interval $[-\pi,\pi)$ whenever convenient.

We denote by $\mathcal{P}(\mathbb{S}^{1})$ the space of periodic probability
measures endowed with the periodic $2$-Wasserstein distance
\begin{equation}
\mathbf{d}_{per}^{2}(\mu,\rho)=\inf\left\{  \int_{\mathbb{S}^{1}%
\times\mathbb{S}^{1}}d_{per}^{2}(x,y)\ d\gamma(x,y):\gamma\in\Gamma(\mu
,\rho)\right\}  , \label{dist-1}%
\end{equation}
where $\Gamma(\mu,\rho)$ stands for the set of probability measures with
marginals $\mu$ and $\rho$, and $d_{per}$ denotes the geodesic distance in
$\mathbb{S}^{1}$. The subspace of absolutely continuous measures in
$\mathcal{P}(\mathbb{S}^{1})$ is denoted by $\mathcal{P}_{ac}(\mathbb{S}%
^{1}).$ From \cite{McCann-2}, for $\mu,\rho\in\mathcal{P}(\mathbb{S}^{1})$
with $\mu\in\mathcal{P}_{ac}(\mathbb{S}^{1})$, there exists an optimal
transport map $\mathbf{t}_{\mu}^{\rho}:\mathbb{S}^{1}\rightarrow\mathbb{S}%
^{1}$ for the Monge problem with quadratic distance cost. Indeed it is
possible to show that $\mathbf{t}_{\mu}^{\rho}$ exists if $\mu$ does not give
mass on points (see the argument below), i.e., $\mu$ has no atoms.

Following \cite{Carrillo-3}, we can consider the application $\mathbf{t}_{\mu
}^{\rho}$ from $[-\pi,\pi]$ to $\mathbb{S}^{1}$. Define the application
$\widetilde{\mathbf{t}_{\mu}^{\rho}}:[-\pi,\pi]\rightarrow\lbrack-\pi,3\pi]$
given by: $\widetilde{\mathbf{t}_{\mu}^{\rho}}(x)$ is the smallest element in
the equivalence class $\mathbf{t}_{\mu}^{\rho}(x)$ such that $\left\vert
\widetilde{\mathbf{t}_{\mu}^{\rho}}(x)-x\right\vert \leq\pi$. In this case,
the geodesic distance of $x$ and $\mathbf{t}_{\mu}^{\rho}(x)$ coincides with
the Euclidean one between $x$ and $\widetilde{\mathbf{t}_{\mu}^{\rho}}(x)$.
So, considering Euclidean quadratic cost, if ${\tilde{\rho}:=\widetilde
{\mathbf{t}_{\mu}^{\rho}}_{\#}}\mu$ then $\widetilde{\mathbf{t}_{\mu}^{\rho}}$
is the optimal transport map between $\mu$ and $\tilde{\rho}$. Thus,
$\widetilde{\mathbf{t}_{\mu}^{\rho}}$ is monotone and
\[
\widetilde{\mathbf{t}_{\mu}^{\rho}}(\pi)=2\pi+\widetilde{\mathbf{t}_{\mu
}^{\rho}}(-\pi)=:2\pi+a.
\]
It follows from monotonicity that $\widetilde{\mathbf{t}_{\mu}^{\rho}}%
([-\pi,\pi])\subset\lbrack a,a+2\pi]$. In short, we can think in optimal
transports as maps from $[-\pi,\pi]$ to $[a,a+2\pi]$.

Let us also comment on another way to see the space $\mathcal{P}%
(\mathbb{S}^{1})$. Recently, in analogy with the construction of the torus as
a quotient space, the authors of \cite{Gangbo} defined an equivalence relation
in $\mathcal{P}_{2}(\mathbb{R})$ by setting
\begin{equation}
\mu\sim\rho\Leftrightarrow\int_{\mathbb{R}}\zeta d\mu=\int_{\mathbb{R}}\zeta
d\rho,\text{ \ }\forall\zeta\in C^{0}(\mathbb{S}^{1}). \label{rel-equi1}%
\end{equation}
Given $\mu\in\mathcal{P}_{2}(\mathbb{R}),$ there exists a $\hat{\mu}$
equivalent to $\mu$ concentrated in $[-\pi,\pi)$. For that, just take the
pushforward of $\mu$ by the map that sends $x\in\mathbb{R}$ to the (unique)
element of $[x]\cap\lbrack-\pi,\pi)$ (here $[x]$ denotes the equivalence class
in $\mathbb{R}/2\pi\mathbb{Z}$). Indeed $\hat{\mu}$ is the unique
representative of $[\mu]$ with these properties. In the sequel, they showed a
relation between the metrics of $\mathcal{P}_{2}(\mathbb{R})$ and
$\mathcal{P}_{2}(\mathbb{R})/\sim$ similar to that between $\mathbb{R}$ and
$\mathbb{S}^{1}$. Indeed, if $\mathbf{d}_{2}$ is the Wasserstein metric in
$\mathcal{P}_{2}(\mathbb{R})$ then
\begin{align}
\mathbf{d}_{per}^{2}(\mu,\rho)  &  =\min\{\mathbf{d}_{2}^{2}(\mu,\rho^{\ast
}):\rho\sim\rho^{\ast}\}\label{relmetric}\\
&  =\min\{\mathbf{d}_{2}^{2}(\mu^{\ast},\rho^{\ast}):\mu\sim\mu^{\ast}\text{
and }\rho\sim\rho^{\ast}\}. \label{realmetric-00}%
\end{align}
We observe that the minimum in (\ref{relmetric}) is reached with the map
$\widetilde{\mathbf{t}_{\mu}^{\rho}}$ (built above), when $\mu$ is supported
in $[-\pi,\pi)$. In fact, one can check that $\tilde{\rho}\sim\rho$ and
\begin{align*}
\mathbf{d}_{2}^{2}(\mu,\tilde{\rho})  &  \leq\int_{\lbrack-\pi,\pi
)}|x-\widetilde{\mathbf{t}_{\mu}^{\rho}}(x)|^{2}\ d\mu(x)\\
&  =\int_{[-\pi,\pi)}d_{per}^{2}(x,t_{\mu}^{\rho}(x)){\ d\mu(x)}\\
&  =\mathbf{d}_{per}^{2}(\mu,\rho),
\end{align*}
as desired. By (\ref{relmetric}), there exist $\gamma\in\mathcal{P}%
_{2}(\mathbb{R}^{2})$ and $\rho^{\ast}\in\mathcal{P}_{2}(\mathbb{R})$ such
that
\begin{equation}
\mathbf{d}_{per}^{2}(\mu,\rho)=\int_{\mathbb{R}^{2}}|x-y|^{2}\ d\gamma
(x,y),\ \gamma\in\Gamma(\mu,\rho^{\ast})\text{ and}\ \rho^{\ast}\sim\rho.
\label{relmetric1}%
\end{equation}
Finally, since $\mathbb{S}^{1}$ is compact, let us remark that the weak
topology (narrow)\textbf{ }in $\mathcal{P}(\mathbb{S}^{1})$ coincides with
that induced by the $p$-Wasserstein metric. In particular, by Prokhorov lemma,
$\mathcal{P}(\mathbb{S}^{1})$ is a compact metric space. The above two ways of
seeing the space $\mathcal{P}(\mathbb{S}^{1})$ will be exploited by us
throughout the paper.

\subsection{\bigskip Geodesics with Constant Velocity in $\mathcal{P}%
(\mathbb{S}^{1})$}

\hspace{0.4cm} In this section, we consider the elements of $\mathbb{S}^{1}$
as equivalence classes $[x]$. In what follows, we recall the definition of a
type of geodesic in metric spaces.

\begin{defn}
Let $(X,d)$ be a metric space. A curve $\varphi:[0,1]\rightarrow X$ is called
a geodesic with constant velocity if
\[
d(\varphi(s),\varphi(t))=|t-s|d(\varphi(0),\varphi(1))\ \ \forall
s,t\in\lbrack0,1].
\]

\end{defn}

Here we present an explicit construction of a geodesic with constant velocity
connecting two arbitrary measures in $\mathcal{P}(\mathbb{S}^{1})$. Define the
multivalued map $g_{t}:\mathbb{S}^{1}\times\mathbb{S}^{1}\rightarrow
\mathbb{S}^{1}$ given by $g_{t}([x],[y])=[(1-t)\hat{x}+t\hat{y}]$, where
\[
(\hat{x},\hat{y})\in\text{Argmin}\{|\bar{x}-\bar{y}|:\bar{x}\in\lbrack
x],\ \bar{y}\in\lbrack y]\}.
\]
Note that $g_{t}$ is a function well defined in the set
\[
\left\{  ([x],[y])\in\mathbb{S}^{1}\times\mathbb{S}^{1}:d_{per}([x],[y])<\pi
\right\}  .
\]
If $d_{per}([x],[y])=\pi$ then we can redefine $g_{t}$ by requiring that
$\hat{y}<\hat{x}$. It follows that
\[
d_{per}(g_{t}([x],[y]),g_{s}([x],[y]))=|t-s|d_{per}([x],[y]),\text{ for all
}s,t\in\lbrack0,1].
\]

Given $\mu_{0},\mu_{1}\in\mathcal{P}(\mathbb{S}^{1})$ and $\gamma\in\Gamma
(\mu_{0},\mu_{1})$ an optimal plan for $\mathbf{d}_{per}(\mu_{0},\mu_{1})$, we define%

\begin{equation}
\mu_{t}:={g_{t}}_{\#}\gamma,\ t\in\lbrack0,1]. \label{geodesic}%
\end{equation}

\begin{prop}
Let $\mu_{0},\mu_{1}\in\mathcal{P}(\mathbb{S}^{1})$. Then, the curve defined
in (\ref{geodesic}) is a geodesic with constant velocity with respect to the
Wasserstein metric in $\mathcal{P}(\mathbb{S}^{1})$.
\end{prop}

\emph{Proof.} Define $\gamma_{t,s}=(g_{t},g_{s})_{\#}\gamma\in\Gamma(\mu
_{t},\mu_{s})$. It follows that
\begin{align}
\mathbf{d}_{per}^{2}(\mu_{t},\mu_{s})  &  \leq\int_{\mathbb{S}^{1}%
\times\mathbb{S}^{1}}{d_{per}^{2}(g_{t}([x],[y]),g_{s}([x],[y]))\ d\gamma
}\label{aux-1}\\
&  =(t-s)^{2}\int_{\mathbb{S}^{1}\times\mathbb{S}^{1}}{d_{per}^{2}%
([x],[y])\ d\gamma},\nonumber
\end{align}
and then $\mathbf{d}_{per}(\mu_{t},\mu_{s})\leq|t-s|\mathbf{d}_{per}(\mu
_{0},\mu_{1})$. If there exist $s,t\in\lbrack0,1]$ such that $s<t$ and the
inequality in (\ref{aux-1}) is strict, then
\begin{align*}
\mathbf{d}_{per}(\mu_{0},\mu_{1})  &  <s\mathbf{d}_{per}(\mu_{0},\mu
_{1})+|t-s|\mathbf{d}_{per}(\mu_{0},\mu_{1})+(1-t)\mathbf{d}_{per}(\mu_{0}%
,\mu_{1})\\
&  =\mathbf{d}_{per}(\mu_{0},\mu_{1}),
\end{align*}
which gives a contradiction. Therefore, we have indeed an equality in
(\ref{aux-1}), as required. \fin

\subsection{Gradient-Flow Structure}

\hspace{0.4cm} Formally, we can write (\ref{eq-1}) and (\ref{eq1-viscous}) as
\begin{align}
u_{t}  &  =\left[  u\left(  \nu\frac{u_{x}}{u}-\mathcal{H}(u)\right)  \right]
_{x}\nonumber\\
&  =\left[  u\left(  \nu\log u-\frac{1}{\pi}\log\left\vert \sin
(x/2)\right\vert \ast u\right)  _{x}\right]  _{x},\text{ for }\nu\geq0.
\label{eq-flow-1}%
\end{align}
Since we are looking for solutions in $\mathcal{P}(\mathbb{S}^{1}),$ equation
(\ref{eq-flow-1}) suggests to define the interaction kernel as
\begin{equation}
W(x)=\left\{
\begin{array}
[c]{cl}%
-\frac{1}{\pi}\log\left\vert \sin(x/2)\right\vert  & \text{ if }x\in
\lbrack-\pi,\pi),\ x\neq0;\\
\infty & \text{ if }x=0;\\
W(x+2\pi)=W(x) & \text{ for }x\in\mathbb{R}.
\end{array}
\right.  \label{interactionpotential}%
\end{equation}
Now consider the free energy functional $\mathcal{F}_{\nu}:\mathcal{P}%
(\mathbb{S}^{1})\rightarrow(-\infty,\infty]$ defined in the following way:
\begin{align}
\mathcal{F}_{\nu}[\mu]  &  =\nu\int_{\lbrack-\pi,\pi)}\mu\log\mu dx+\int
\int_{[-\pi,\pi)^{2}}W(x-y)d\mu(x)d\mu(y)\label{energyfunctional}\\
&  =:\nu\mathcal{U}[\mu]+\mathcal{F}_{0}[\mu],\text{ for all }\mu
\in\mathcal{P}_{ac}(\mathbb{S}^{1})\text{ and }\nu>0.\nonumber
\end{align}
Here we are identifying an absolutely continuous measure with your density
with respect to Lebesgue measure. For $\mu\in\mathcal{P}(\mathbb{S}%
^{1})\backslash\mathcal{P}_{ac}(\mathbb{S}^{1})$ and $\nu>0$, $\mathcal{F}%
_{\nu}[\mu]=\mathcal{U}[\mu]=\infty.$ In the case $\nu=0$ we simply define the
functional as
\begin{equation}
\mathcal{F}_{\nu=0}[\mu]=\mathcal{F}_{0}[\mu],\text{ for all }\mu
\in\mathcal{P}(\mathbb{S}^{1}). \label{energyfunctional1}%
\end{equation}
We recall that the domain of the functional, denoted by $D(\mathcal{F_{\nu}}%
)$, is defined as the set
\begin{equation}
D(\mathcal{F}_{\nu})=\left\{  \ \mu\in\mathcal{P}(\mathbb{S}^{1}%
);\mathcal{F}_{\nu}(\mu)<\infty\text{ }\right\}  . \label{dom-1}%
\end{equation}
So, we can write (\ref{eq-flow-1}) in the form (\ref{continuity}) with
$\upsilon=\frac{\delta\mathcal{F}_{\nu}}{\delta u},$ that is
\begin{equation}
u_{t}=\left[  u\left(  \frac{\delta\mathcal{F}_{\nu}}{\delta u}\right)
_{x}\right]  _{x} \label{grad-flow-1}%
\end{equation}
which has the structure of gradient flow in $\mathcal{P}(\mathbb{S}^{1})$
corresponding to the energy functional $\mathcal{F}_{\nu}$.

\begin{rem}
\label{welldefinition1}Let us remark that if we adopt the interpretation of
$\mathcal{P}(\mathbb{S}^{1})$ as $\mathcal{P}_{2}(\mathbb{R})/\sim,$ the
definition of $\mathcal{F}_{\nu}$ is as follows:

\begin{itemize}
\item[(i)] For the entropy part of $\mathcal{F}_{\nu}$ and $\mu\in
\mathcal{P}_{2}(\mathbb{R})$, we consider the unique $\mu^{\ast}\sim\mu$ such
that $\mu^{\ast}$ is supported in $[-\pi,\pi)$ and define $\mathcal{U}%
[\mu]:=\mathcal{U}[\mu^{\ast}]$ when $\mu\in\mathcal{P}_{ac}(\mathbb{S}^{1})$,
and $\mathcal{U}[\mu]=\infty$ otherwise.

\item[(ii)] For the interaction part, first observe that
\[
\int\!\int_{\mathbb{R}^{2}}W(x-y)\ d\mu(x)d\mu(y)=\int\!\int_{\mathbb{R}^{2}%
}W(x-y)\ d\mu^{\ast}(x)d\mu^{\ast}(y),
\]
for any $\mu\sim\mu^{\ast}$ in $\mathcal{P}_{2}(\mathbb{R})$. In fact, this
equality follows by approximating the kernel $W$ monotonically from below by
periodic continuous functions and then applying the monotone convergence
theorem. Therefore, we define the interaction functional by
\[
\mathcal{F}_{0}[\mu]=\int\!\int_{\mathbb{R}^{2}}W(x-y)\ d\mu^{\ast}%
(x)d\mu^{\ast}(y).
\]

\end{itemize}
\end{rem}

Connected to item (i) in the previous remark, let us show some kind of
invariance for the entropy functional $\mathcal{U}[\mu]$.

\begin{lem}
\label{entropyinvariance} Let $U:[0,\infty)\rightarrow\mathbb{R}$ with
$U(0)=0$, $\mu,\rho\in\mathcal{P}_{ac}(\mathbb{R})$ such that $\mu\sim\rho$,
$\mu$ is supported in $[a,a+2\pi)$ and $\rho$ is supported in $[b,b+2\pi)$.
Then, if $f$ and $h$ are the densities of $\mu$ and $\rho$ respectively, we
have $U\circ f\in L^{1}(\mathbb{R},dx)$ if and only if $U\circ h\in
L^{1}(\mathbb{R},dx)$ and
\[
\int_{\lbrack a,a+2\pi)}{U\circ f\ dx}=\int_{[b,b+2\pi)}{U\circ h\ dx}.
\]

\end{lem}

\emph{Proof.} Without loss of generality, we can assume that $\mu$ is
supported in $[0,2\pi)$. Define the function $T:[b,b+2\pi)\rightarrow
\lbrack0,2\pi)$ by
\[
T(x)=\left\{
\begin{array}
[c]{ll}%
x-2\pi\lfloor\frac{b}{2\pi}\rfloor & \text{ if }b\leq x<2\pi\left(
1+\lfloor\frac{b}{2\pi}\rfloor\right)  ;\\
x-2\pi\left(  1+\lfloor\frac{b}{2\pi}\rfloor\right)  & \text{ if }2\pi\left(
1+\lfloor\frac{b}{2\pi}\rfloor\right)  \leq x<b+2\pi,
\end{array}
\right.
\]
where $\lfloor\cdot\rfloor$ stands for the greatest integer function. Then,
clearly $T$ is bijective and it is straightforward to check that $T_{\#}%
\rho\sim\rho$ and $T_{\#}\rho$ is supported in $[0,2\pi)$. Since there exists
a unique representative equivalent to $\rho$ supported in $[0,2\pi)$ (see
Section 2.1), we conclude that $T_{\#}\rho=\mu$. Let $f$ and $h$ be the
densities of $\mu$ and $\rho$, respectively. Then, for any $\zeta\in
C(\mathbb{S}^{1}),$ we have
\[
\int_{\lbrack0,2\pi)}\zeta(x)h(T^{-1}(x))\ dx=\int_{[0,2\pi)}\zeta
(x)f(x)\ dx.
\]
It follows that $h\circ T^{-1}=f$ a.e. in $[0,2\pi)$ and a change of variables
completes the proof.\fin

\begin{rem}
\label{domW} A natural question is to know how much large in $\mathcal{P}%
(\mathbb{S}^{1})$ the set (\ref{dom-1}) is when $\nu=0$. In fact,
$\mathcal{F}_{0}[\mu]=+\infty$ when $\mu\in\mathcal{P}(\mathbb{S}%
^{1})\backslash\mathcal{P}_{ac}(\mathbb{S}^{1})$ gives mass on points. On the
other hand, it is interesting to note that $\mathcal{F}_{0}$ \ may be finite
on singular measures (with respect to Lebesgue one) that concentrate mass on
sets with positive Hausdorff dimension (see proposition below).
\end{rem}

\bigskip

\begin{prop}
\label{dom-set} For $s=\log(2)/\log(3),$ let $\mathcal{H}^{s}$ denote the
$s$-dimensional Hausdorff measure in $\mathbb{R}$ and let $\mathbf{C}$ denote
the Cantor ternary set in $[-\pi,\pi]$. Let $\mu^{s}\in\mathcal{P}%
(\mathbb{S}^{1})$ be defined by
\begin{equation}
\mu^{s}(A)=\frac{1}{2\pi}\mathcal{H}^{s}(A\cap\mathbf{C}), \label{measure}%
\end{equation}
for all Lebesgue measurable set $A$. \ Then, $\mu^{s}$ is singular with
respect to the Lebesgue measure and $\mathcal{F}_{0}[\mu^{s}]<\infty$.
\end{prop}

\emph{Proof.} The Hausdorff dimension of $\mathbf{C}$ is $s:=\log(2)/\log(3)$
and $\mathcal{H}^{s}(\mathbf{C})=2\pi$ (see \cite[p. 34]{Falconer}). Define
the probability measure $\mu^{s}(A)=\frac{1}{2\pi}\mathcal{H}^{s}%
(A\cap\mathbf{C})$ in $\mathbb{S}^{1}\equiv\lbrack-\pi,\pi)$. Now notice that
\[
-\log\left\vert \sin\left(  \frac{x}{2}\right)  \right\vert \leq
-\log\left\vert \frac{x}{4}\right\vert ,
\]
for $|x|\leq x_{0}$ and some $x_{0}>0$. Thus
\begin{align}
\mathcal{F}_{0}[\mu^{s}]  &  \leq\frac{-1}{\pi}\int\!\int_{[|x-y|\leq
x_{0}]\cap\lbrack-\pi,\pi)^{2}}\log\left\vert \frac{x-y}{4}\right\vert
d\mu^{s}(y)d\mu^{s}(x)\nonumber\\
&  +\frac{-1}{\pi}\int\!\int_{[|x-y|>x_{0}]\cap\lbrack-\pi,\pi)^{2}}%
\log\left\vert \sin(\frac{x-y}{2})\right\vert d\mu^{s}(y)d\mu^{s}(x).
\label{aux-int-2}%
\end{align}
The second integral in (\ref{aux-int-2}) is obviously finite, while the first
can be bounded after showing%
\begin{equation}
\int_{\lbrack-\pi,\pi)}-\log\left\vert \frac{x-y}{2\pi}\right\vert d\mu(y)\leq
C, \label{aux-int-3}%
\end{equation}
for some universal constant $C>0$ . The estimate (\ref{aux-int-3}) can be
showed by approximating monotonically from below the function $-\log\left\vert
\frac{x-\cdot}{2\pi}\right\vert $ by functions based on the construction of
the Cantor set and then by applying the monotone convergence theorem.\fin

\begin{rem}
\label{dom-link-trans}There is an optimal transport map $\mathbf{t}_{\mu
}^{\rho}$ between $\mu$ and $\rho,$ when $\mu$ does not give mass to points
(see \cite[p.75]{Villani}). Then, the subject of Remark \ref{domW} is not a
restriction to the existence of an optimal transport map between elements of
$D(\mathcal{F_{\nu}}).$
\end{rem}

\bigskip

\section{Properties of the Functionals}

\hspace{0.4cm} The aim of this section is to obtain lower semicontinuity,
coercivity and convexity properties for the free energy functional
$\mathcal{F}_{\nu}$.

\subsection{Lower semicontinuity}

\hspace{0.4cm} In the next lemma we adapt some ideas of \cite{McCann} for our context.

\begin{lem}
\label{semicontinuidade} The functional $\mathcal{F}_{\nu}$ is lower
semicontinuous with respect to the weak topology that coincides with the
topology induced by the Wasserstein metric.
\end{lem}

\emph{Proof.} In view of (\ref{energyfunctional}), it is sufficient to prove
the lower semicontinuity of $\mathcal{U}$ and $\mathcal{F}_{0}$.

\textit{Step 1 }(Semicontinuity of $\mathcal{U}$): Denote by $U(t)=t\log
t,\ t\geq0$. We start with the lower semicontinuity of the functional
$\mathcal{U}$. Let $\mu_{k}\rightarrow\mu$ be a sequence weakly convergent in
$\mathcal{P}(\mathbb{S}^{1})$ and assume that $\mu\in\mathcal{P}%
_{ac}(\mathbb{S}^{1})$. We can assume that $\liminf_{k\rightarrow\infty
}\mathcal{U}(\mu_{k})<\infty$ and then $\mu_{k}\in\mathcal{P}_{ac}%
(\mathbb{S}^{1})$ with $d\mu_{k}=\rho_{k}dx$. Let $\eta$ be a nonnegative
smooth compactly supported mollifier with mass equal to $1$ in $\mathbb{S}%
^{1}$. Let $\eta_{\delta}(x)=\frac{1}{\delta}\eta(\frac{x}{\delta})$, for
$\delta>0$. Then, by Jensen inequality, we have
\begin{equation}
\int_{-\pi}^{\pi}\rho_{k}(y)\log\rho_{k}(y)\eta_{\delta}(x-y)dy\geq U\left(
\int_{-\pi}^{\pi}\rho_{k}(y)\eta_{\delta}(x-y)dy\right)  . \label{sci1}%
\end{equation}
It follows by integrating in $x$ that
\begin{align}
\liminf_{k\rightarrow\infty}\mathcal{U}(\mu_{k})  &  \geq\liminf
_{k\rightarrow\infty}\mathcal{U}\left(  \rho_{k}\ast\eta_{\delta}\right)
\label{sci2}\\
&  \geq\int_{-\pi}^{\pi}\liminf_{k\rightarrow\infty}U\left(  \rho_{k}\ast
\eta_{\delta}(x)\right)  dx\nonumber\\
&  =\int_{-\pi}^{\pi}U\left(  \rho\ast\eta_{\delta}(x)\right)  dx,\nonumber
\end{align}
where we have used Fatou lemma in the second inequality. The equality follows
from the weak convergence of $\mu_{k}$ and continuity of the function $U$.
Since the support of $\eta_{\delta}$ is shrinking to a point, as
$\delta\rightarrow0$, then we obtain from Lebesgue differentiation theorem
that
\[
\lim_{\delta\rightarrow0}\rho\ast\eta_{\delta}(x)=\rho(x)\text{ a.e. in
}\mathbb{S}^{1}.
\]
The semicontinuity follows by applying again Fatou lemma and continuity of
$U$. Let us assume now that $\mathcal{U}(\mu)=\infty$. Without loss of
generality, we can assume that $\mathcal{U}(\mu_{k})<\infty,$ for $k$ large
enough. Let $d\mu=d\mu_{sing}+\rho dx$ a Lebesgue decomposition in singular
and absolutely continuous parts. By regularity of $\mu_{sing},$ there exist a
compact $K\subset\mathbb{S}^{1}$ with null Lebesgue measure and a number $m>0$
such that $\mu_{sing}(K)\geq m>0$. Thus, there exists an open set $O$
containing $K$ with Lebesgue measure arbitrarily small. Recall that weak
convergence implies $\liminf_{k\rightarrow\infty}\mu_{k}(O)\geq\mu(O)$. We
denote by $|O|$ the Lebesgue measure of $O$ and $c_{0}=\inf_{t\geq0}U(t)$. By
Jensen inequality, we have
\begin{align*}
\mathcal{U}[\mu_{k}]-2\pi c_{0}  &  =\int_{-\pi}^{\pi}{(}U(\rho_{k}%
(x))-c_{0}{)dx}\\
&  \geq|O|\left[  U\left(  \int_{O}\frac{\rho_{k}}{|O|}dx\right)
-c_{0}\right] \\
&  =|O|\left[  U\left(  \frac{\mu_{k}(O)}{|O|}\right)  -c_{0}\right]  .
\end{align*}
Since $U(t)/t$ is increasing and $U(t)/t\rightarrow+\infty$ as $t\rightarrow
\infty$, we obtain
\begin{align*}
\mathcal{U}[\mu_{k}]-2\pi c_{0}  &  \geq\mu_{k}(O)\frac{|O|}{m}U\left(
\frac{m}{|O|}\right)  -c_{0}|O|\\
&  \geq|O|U\left(  \frac{m}{|O|}\right)  -c_{0}|O|,
\end{align*}
for $k$ large enough. We conclude by letting $|O|\rightarrow0$.

\textit{Step 2 }(Semicontinuity of $\mathcal{F}_{0}$): Note that $W$ can be
approximated monotonically from below by periodic bounded continuous functions
$W_{l}$. Thus, if $\mu_{k}\rightarrow\mu$ weakly then
\begin{align*}
\int\!\int_{[-\pi,\pi)^{2}}W_{l}(y-x)d\mu(y)d\mu(x)  &  =\lim_{k\rightarrow
\infty}\int\!\int_{[-\pi,\pi)^{2}}W_{l}(y-x)d\mu_{k}(y)d\mu_{k}(x)\\
&  \leq\liminf_{k\rightarrow\infty}\int\!\int_{[-\pi,\pi)^{2}}W(y-x)d\mu
_{k}(y)d\mu_{k}(x),
\end{align*}
because the weak convergence of $\mu_{k}$ implies weak convergence of $\mu
_{k}\times\mu_{k}\rightarrow\mu\times\mu$. Now, an application of the monotone
convergence theorem finishes the proof. \fin

\subsection{Existence of Minimizer}

\hspace{0.4cm} At this level, we can show the existence of a minimizer for
$\mathcal{F}_{\nu}$. In fact, it is bounded from below, because
\[
\nu\mathcal{U}[\mu]\geq2\pi\nu\inf_{t>0}t\log(t)=-2\pi\nu e^{-1}%
\]
and
\[
\mathcal{W}[\mu]\geq0.
\]
Choose a minimizer sequence $\mu_{k}$ for $\inf\mathcal{F}_{\nu}$. In view of
the weak compactness of $\mathcal{P}(\mathbb{S}^{1}),$ we can assume (up to a
subsequence) that $\mu_{k}\rightarrow\mu_{0}$ in $\mathcal{P}(\mathbb{S}^{1}%
)$. It follows from the lower semicontinuity that
\[
\inf\mathcal{F}_{\nu}\leq\mathcal{F}_{\nu}[\mu_{0}]\leq\liminf_{k\rightarrow
\infty}\mathcal{F}_{\nu}[\mu_{k}]=\inf\mathcal{F}_{\nu},
\]
and so $\mu_{0}$ is a minimizer of $\mathcal{F}_{\nu}$.

\subsection{Convexity in $\mathcal{P}(\mathbb{S}^{1})$}

\hspace{0.4cm} The minimizer obtained in the previous section is unique since
we show some kind of convexity for $\mathcal{F}_{\nu}$.

Given $\mu_{0},\mu_{1},\omega\in\mathcal{P}(\mathbb{S}^{1})$, we know that
there exist $\mu_{0}^{\ast}\sim\mu_{0},$ $\mu_{1}^{\ast}\sim\mu_{1}$ and plans
$\gamma_{0}\in\Gamma(\omega,\mu_{0}^{\ast})$, $\gamma_{1}\in\Gamma(\omega
,\mu_{1}^{\ast})$ such that (\ref{relmetric1}) is true. In the case when
$\omega$ is supported in $[-\pi,\pi)$ and does not give mass to points, we
know that it is possible to choose $\mu_{i}^{\ast}=\widetilde{\mathbf{t}%
_{\omega}^{\mu_{i}}}_{\#}\omega$ and $\gamma_{i}=(I,\widetilde{\mathbf{t}%
_{\omega}^{\mu_{i}}})_{\#}\omega$,\ $i=0,1$, where $\widetilde{\mathbf{t}%
_{\omega}^{\mu_{i}}}$ is the map built in Section 2.1. We can consider it as
$\widetilde{\mathbf{t}_{\omega}^{\mu_{i}}}:[-\pi,\pi)\rightarrow\lbrack
a_{i},a_{i}+2\pi)$, $i=0,1$. Thus, the map
\[
(1-t)\mathbf{t}_{\omega}^{\mu_{0}}+t\mathbf{t}_{\omega}^{\mu_{1}}:[-\pi
,\pi)\rightarrow\lbrack(1-t)a_{0}+ta_{1},(1-t)a_{0}+ta_{1}+2\pi)
\]
can be seen as a map from $\mathbb{S}^{1}$ to $\mathbb{S}^{1}.$

In what follows, following the terminology of \cite{Ambrosio}, we define the
concepts of generalized geodesic and convexity along generalized geodesics in
$\mathcal{P}(\mathbb{S}^{1})$.

\begin{defn}
\label{Def-Geo-1}Given $\mu_{0}$, $\mu_{1}$, $\omega\in\mathcal{P}%
(\mathbb{S}^{1})$, choose pairs $(\mu_{0}^{\ast},\gamma_{0})$ and $(\mu
_{1}^{\ast},\gamma_{1})$ such that $\gamma_{0}\in\Gamma(\omega,\mu_{0}^{\ast
})$ and $\gamma_{1}\in\Gamma(\omega,\mu_{1}^{\ast})$ and (\ref{relmetric1}) is
valid. A generalized geodesic connecting $\mu_{0}$ to $\mu_{1}$ in
$\mathcal{P}(\mathbb{S}^{1})$, with base point in $\omega$ and induced by
$\mbox{\boldmath{$\gamma$}}$, is a curve of equivalence classes of the type
$\mu_{t}^{g}:=((1-t)P_{2}+tP_{3})_{\#}\mbox{\boldmath{$\gamma$}}$,
$t\in\lbrack0,1],$ where $P_{2}$ and $P_{3}$ are the second and third
projections from $\mathbb{R}^{3}$ onto $\mathbb{R}$, and
$\mbox{\boldmath{$\gamma$}}\in\Gamma(\omega,\mu_{0}^{\ast},\mu_{1}^{\ast})$ is
such that $(P_{1},P_{2})_{\#}\mbox{\boldmath{$\gamma$}}=\gamma_{0}$ and
$(P_{1},P_{3})_{\#}\mbox{\boldmath{$\gamma$}}=\gamma_{1}$.
\end{defn}

\begin{rem}
\label{rem-0} Notice that the convex combination $(1-t)P_{2}+tP_{3}$ can be
outside of $\mathbb{S}^{1}$. Then, $\mu_{t}^{g}$ should be understood as a
curve of equivalence classes, according to the identification between
$\mathcal{P}(\mathbb{S}^{1})$ and $\mathcal{P}_{2}(\mathbb{R})/\sim$.
\end{rem}

\begin{rem}
\label{rem-1}When $\omega$, $\mu_{0}$ and $\mu_{1}$ are supported in
$[-\pi,\pi)$ and $\omega$ does not give mass to points, we can take
$\mbox{\boldmath{$\gamma$}}=(I,\widetilde{\mathbf{t}_{\omega}^{\mu_{0}}%
},\widetilde{\mathbf{t}_{\omega}^{\mu_{1}}})_{\#}\omega$. Therefore, a
generalized geodesic in $\mathcal{P}(\mathbb{S}^{1})$ is given by the
equivalence class of $\mu_{t}^{g}=((1-t)\widetilde{\mathbf{t}_{\omega}%
^{\mu_{0}}}+t\widetilde{\mathbf{t}_{\omega}^{\mu_{1}}})_{\#}\omega$.
\end{rem}

\begin{defn}
\label{convexity}We say that a functional $\mathcal{F}:\mathcal{P}%
(\mathbb{S}^{1})\rightarrow(-\infty,\infty]$ is $\lambda-$convex along of
generalized geodesics, for some $\lambda\in\mathbb{R}$, if given $\omega
,\mu_{0},\mu_{1}\in D(\mathcal{F})$ (the domain of $\mathcal{F}$) there exists
a generalized geodesic $\mu_{t}^{g}$ connecting $\mu_{0}$ to $\mu_{1}$, based
in $\omega$ and induced by $\mbox{\boldmath{$\gamma$}}$, such that
\[
\mathcal{F}[\mu_{t}^{g}]\leq(1-t)\mathcal{F}[\mu_{0}]+t\mathcal{F}[\mu
_{1}]-\frac{\lambda}{2}t(1-t)\mathbf{d}_{\mbox{\boldmath{$\gamma$}}}^{2}%
(\mu_{0}^{\ast},\mu_{1}^{\ast}),
\]
where $\mu_{i}^{\ast}\sim\mu_{i}$ and $\mathbf{d}_{\mbox{\boldmath{$\gamma$}}}%
^{2}(\mu_{0}^{\ast},\mu_{1}^{\ast})=\int_{\mathbb{R}^{3}}|x_{2}-x_{3}%
|^{2}\ d\mbox{\boldmath{$\gamma$}}(x_{1},x_{2},x_{3}).$
\end{defn}

\bigskip

\begin{rem}
\label{rem-3-aux}In Definition \ref{convexity}, from (\ref{realmetric-00})
notice that $\mathbf{d}_{\mbox{\boldmath{$\gamma$}}}^{2}(\mu_{0}^{\ast}%
,\mu_{1}^{\ast})\geq\mathbf{d}_{2}^{2}(\mu_{0}^{\ast},\mu_{1}^{\ast}%
)\geq\mathbf{d}_{per}^{2}(\mu_{0},\mu_{1}),$ for all $\mu_{0},\mu_{1}%
\in\mathcal{P}(\mathbb{S}^{1}).$
\end{rem}

\bigskip

The next lemma shows that $\mathcal{F}_{\nu}$ is convex along generalized geodesics.

\begin{lem}
\label{convexityfunctional} The energy functional $\mathcal{F}_{\nu}$ is
strictly $0-$convex along generalized geodesics. Thus, the minimizer of
$\mathcal{F}_{\nu}$ obtained in Section 3.2 is unique.
\end{lem}

\emph{Proof.} Let $\omega,\mu_{0},\mu_{1}\in D(\mathcal{F}_{\nu})$ be
supported in $[-\pi,\pi)$. Since $\omega,\mu_{0}$ and $\mu_{1}$ have no atoms,
we can choose the generalized geodesic with representative $\mu_{t}%
^{g}=((1-t)\widetilde{\mathbf{t}_{\omega}^{\mu_{0}}}+t\widetilde
{\mathbf{t}_{\omega}^{\mu_{1}}})_{\#}\omega$ (see Remark \ref{rem-1}).

For the entropy functional $\mathcal{U}$, we know that $\omega,\mu_{0},\mu
_{1}\in\mathcal{P}_{ac}(\mathbb{S}^{1})$ and have that $\mu_{t}^{g}%
\in\mathcal{P}_{ac}(\mathbb{S}^{1})$ with support contained in $[(1-t)a_{0}%
+ta_{1},(1-t)a_{0}+ta_{1}+2\pi)$. Denoting by $f_{t}$ the density of $\mu
_{t}^{g}$ and using Lemma \ref{entropyinvariance}, we have that the entropy is
given by
\begin{equation}
\mathcal{U}[\mu_{t}^{g}]=\int_{\mathbb{R}}f_{t}\log(f_{t})\ dx.
\label{aux-func-1}%
\end{equation}
Now the convexity follows by arguing as in \cite[Proposition 9.3.9]{Ambrosio}
(see also \cite{McCann} for the case of displacement interpolation curves).

On the other hand, recall that $W(x)=-\frac{1}{\pi}\log\left\vert \sin\left(
\frac{x}{2}\right)  \right\vert $ (for$\ x\neq0$) is the kernel of
$\mathcal{F}_{0}$ and we have
\[
W^{\prime\prime}(x)=\frac{1}{4\pi}\csc^{2}(\frac{x}{2})>0,
\]
for $-2\pi<x<2\pi$ with $x\neq0.$ Thus, $W$ is convex on the segments
$(-2\pi,0)$ and $(0,2\pi)$. Now, a modification of an argument in
\cite{Car-Fer-Pre} can be used in order to show that $\mathcal{F}_{0}$ is
convex along generalized geodesics. In fact, by Remark \ref{welldefinition1}
(ii), we have that
\[
\mathcal{F}_{0}[\mu_{t}^{g}]=\int\!\int_{[-\pi,\pi)^{2}}{W((1-t)(\widetilde
{\mathbf{t}_{\omega}^{\mu_{0}}}(x)-\widetilde{\mathbf{t}_{\omega}^{\mu_{0}}%
}(y))+t(\widetilde{\mathbf{t}_{\omega}^{\mu_{1}}}(x)-\widetilde{\mathbf{t}%
_{\omega}^{\mu_{1}}}(y)))\ d(\omega\times\omega).}%
\]
Note that $(x,y)\in\lbrack-\pi,\pi)^{2}$ implies
\[
-2\pi<\widetilde{\mathbf{t}_{\omega}^{\mu_{i}}}(x)-\widetilde{\mathbf{t}%
_{\omega}^{\mu_{i}}}(y)<2\pi,\ i=0,1,
\]
and monotonicity of Euclidean transport maps implies that $\widetilde
{\mathbf{t}_{\omega}^{\mu_{0}}}(x)-\widetilde{\mathbf{t}_{\omega}^{\mu_{0}}%
}(y)\geq0$ if and only if $\widetilde{\mathbf{t}_{\omega}^{\mu_{1}}%
}(x)-\widetilde{\mathbf{t}_{\omega}^{\mu_{1}}}(y)\geq0$. Using the convexity
of $W$ on $(-2\pi,0)$ and $(0,2\pi)$ separately, we are done. \fin

\begin{rem}
\label{convexityongeodesics}For the functionals above, it is possible to show
convexity along geodesics with constant velocity (see Section 2.2) instead of
generalized geodesics.
\end{rem}

\bigskip

The next lemma contains an essential property for the convergence of the Euler scheme.

\begin{lem}
\label{convexitywasserstein} For each fixed $\rho\in\mathcal{P}(\mathbb{S}%
^{1})$, the functional $\rho\rightarrow\mathbf{d}_{per}^{2}(\mu,\rho)$ is
$2$-convex along generalized geodesics.
\end{lem}

\emph{Proof.} Let $\omega,\mu_{0},\mu_{0}^{\ast},\mu_{1},\mu_{1}^{\ast}$ and
$\mbox{\boldmath{$\gamma$}}$ be as in Definition 3.2 such that $\mu_{0}^{\ast
}$ and $\mu_{1}^{\ast}$ are minimum points in (\ref{relmetric}) for
$\mathbf{d}_{per}^{2}(\omega,\mu_{0})$ and $\mathbf{d}_{per}^{2}(\omega
,\mu_{1})$, respectively. Using the $2$-convexity of the $2$-Wasserstein
metric in $\mathcal{P}_{2}(\mathbb{R})$ and (\ref{relmetric}), we obtain
\begin{align}
\mathbf{d}_{per}^{2}(\omega,\mu_{t}^{g})  &  \leq\mathbf{d}_{2}^{2}(\omega
,\mu_{t}^{g})\nonumber\\
&  =(1-t)\mathbf{d}_{2}^{2}(\omega,\mu_{0}^{\ast})+t\mathbf{d}_{2}^{2}%
(\omega,\mu_{1}^{\ast})-t(1-t)\mathbf{d}_{\mbox{\boldmath{$\gamma$}}}^{2}%
(\mu_{0}^{\ast},\mu_{1}^{\ast})\nonumber\\
&  =(1-t)\mathbf{d}_{per}^{2}(\omega,\mu_{0})+t\mathbf{d}_{per}^{2}(\omega
,\mu_{1})-t(1-t)\mathbf{d}_{\mbox{\boldmath{$\gamma$}}}^{2}(\mu_{0}^{\ast}%
,\mu_{1}^{\ast}). \label{aux-22-gen}%
\end{align}
as desired. \fin

\section{Global well-posedness and inviscid limit}

\hspace{0.4cm} In the previous sections we have obtained key properties for
the functional $\mathcal{F}_{\nu}$ in order to construct gradient-flow
solutions via an abstract Euler scheme. In fact, this scheme can be carried
out in general metric spaces since the corresponding functional satisfies
certain conditions (see \cite[Chapter 4]{Ambrosio}).

\subsection{Gradient-Flow Solutions}

\hspace{0.4cm} Here we consider the Euler discrete approximation scheme for
gradient flows in $\mathcal{P}(\mathbb{S}^{1})$. Let $\mu\in\mathcal{P}%
(\mathbb{S}^{1})$ and let $\tau>0$ be the time step. Consider the functional
$\Psi_{\nu}(\tau,\mu;\cdot):\mathcal{P}(\mathbb{S}^{1})\rightarrow
(-\infty,\infty]$ as
\begin{equation}
\Psi_{\nu}(\tau,\mu;\rho):=\frac{1}{2\tau}\mathbf{d}_{per}^{2}(\mu
,\rho)+\mathcal{F}_{\nu}[\rho], \label{euler1}%
\end{equation}
where $\mathbf{d}_{per}$ stands for the $2$-Wasserstein distance in
$\mathcal{P}(\mathbb{S}^{1})$ (see (\ref{dist-1})) and $\mathcal{F}_{\nu}$ is
the functional associated to (\ref{eq1-viscous}), for $\nu\geq0.$

For $\mu_{0}\in\mathcal{P}(\mathbb{S}^{1})$, define the interactive sequence
$(\mu_{\nu,\tau}^{k})_{k=0}^{\infty}$ by
\begin{align}
\mu_{\nu,\tau}^{0}  &  =\mu_{0};\nonumber\\
\mu_{\nu,\tau}^{k}  &  =\text{Argmin}\Psi_{\nu}(\tau,\mu_{\nu,\tau}%
^{k-1};\cdot). \label{scheme-1}%
\end{align}
Again, the minimizer above exists by the compactness of $\mathcal{P}%
(\mathbb{S}^{1})$ and is unique by convexity properties of $\mathcal{F}_{\nu}$
and $\mathbf{d}_{per}^{2}$ (see Section 3.3). This allows to define the
approximate discrete solution of the gradient flow equation (see
\cite{Jordan-Otto} and \cite{Ambrosio})
\begin{align*}
\frac{\partial\mu}{\partial t}  &  =-\text{grad}_{\mathbf{d}_{per}%
}(\mathcal{F}_{\nu})\\
&  =\partial_{x}\cdot(\mu\partial_{x}(\frac{\delta\mathcal{F}_{\nu}}{\delta
\mu}))
\end{align*}
by setting
\begin{equation}
\mu_{\nu,\tau}(t)=\mu_{\nu,\tau}^{k}\text{ if }t\in\lbrack k\tau,(k+1)\tau).
\label{discrete-sol-1}%
\end{equation}

In the following we obtain the well-posedness in $\mathcal{P}(\mathbb{S}^{1})$
and some properties of the gradient flow of $\mathcal{F}_{\nu}.$

\begin{teor}
\label{fluxograd} Let $\nu\geq0$, $\mu_{0}\in\mathcal{P}(\mathbb{S}^{1})$ and
$\mathcal{F}_{\nu}$ be the functional defined in (\ref{energyfunctional}%
)-(\ref{energyfunctional1}).

\begin{itemize}
\item[(i)] The discrete solution $\mu_{\tau}(t)$ defined in
(\ref{discrete-sol-1}) converges locally uniformly to a locally Lipschitz
curve $\mu(t)=S[\mu_{0}](t)$ in $\mathcal{P}(\mathbb{S}^{1})$ which is the
unique gradient flow of $\mathcal{F}_{\nu}$ with $\mu(0)=\mu_{0}$.

\item[(ii)] The map $t\rightarrow S[\mu_{0}](t)$ is a $0$-contracting
semigroup in $\mathcal{P}(\mathbb{S}^{1})$, i.e.
\[
\mathbf{d}_{per}(S[\mu_{0}](t),S[\rho_{0}](t))\leq\mathbf{d}_{per}(\mu
_{0},\rho_{0}),\text{ for }\mu_{0},\rho_{0}\in\mathcal{P}(\mathbb{S}^{1}).
\]

\item[(iii)] If $\nu>0$ then $\mu(t)\in\mathcal{P}_{ac}(\mathbb{S}^{1}),$ for
all $t>0.$ In the case $\nu=0$, the measure solution $u(\cdot,t)=\mu(t)$
concentrates mass at most on sets of positive Hausdorff dimension.

\item[(iv)] Let $\bar{\mu}$ be the unique minimum of $\mathcal{F}_{\nu}$.
Then, the map $t\rightarrow\mathbf{d}_{per}(\mu(t),\bar{\mu})$ is not
increasing and
\begin{equation}
\mathcal{F}_{\nu}(\mu(t))-\mathcal{F}_{\nu}(\bar{\mu})\leq\frac{\mathbf{d}%
_{per}^{2}(\mu_{0},\bar{\mu})}{2t},\text{ for all }t>0. \label{func-conv}%
\end{equation}

\item[(v)] The minimum $\bar{\mu}$ is a stationary gradient flow solution,
i.e., $\bar{\mu}\equiv S[\bar{\mu}](t).$ Moreover, for all $\mu_{0}%
\in\mathcal{P}(\mathbb{S}^{1}),$ $\mu(t)\rightarrow\bar{\mu}$ in
$\mathcal{P}(\mathbb{S}^{1}),$ as $t\rightarrow\infty.$

\item[(vi)] If $\mu_{0}\in D(\mathcal{F}_{\nu})$, then we have the following
error estimate
\begin{equation}
\mathbf{d}_{per}^{2}(\mu_{\tau}(t),\mu(t))\leq\tau\left(  \mathcal{F}_{\nu
}(\mu_{0})+2\pi\nu e^{-1}\right)  . \label{error}%
\end{equation}

\end{itemize}
\end{teor}

\emph{Proof.} We have showed that the energy functional $\mathcal{F}_{\nu}$ is
lower bounded (and so coercive), lower semicontinuous (Lemma
\ref{semicontinuidade}) and $0$-convex along generalized geodesics (Lemma
\ref{convexityfunctional}). On the other hand, Lemmas
\ref{convexityfunctional} and \ref{convexitywasserstein} imply that the
functional defined in (\ref{euler1}) is convex along generalized geodesics.
Finally we recall that $L^{\infty}(\mathbb{S}^{1})\cap\mathcal{P}%
(\mathbb{S}^{1})\subset D(\mathcal{F}_{\nu})$ is dense in $\mathcal{P}%
(\mathbb{S}^{1})$, and therefore we can take $\mu_{0}\in\mathcal{P}%
(\mathbb{S}^{1})$. Now, items (i), (ii),(iv) and (vi) follow from the abstract
theory of \cite[Theorems 4.0.4 and 4.0.7]{Ambrosio} in general metric spaces.
For (iii), we know that $S[\mu_{0}](t)\in D(\mathcal{F}_{\nu})\subset
\mathcal{P}_{ac}(\mathbb{S}^{1})$ for $\nu>0.$ In the case $\nu=0,$ we have
that $S[\mu_{0}](t)\in D(\mathcal{F}_{0})$ and, by Proposition \ref{dom-set},
$D(\mathcal{F}_{0})$ contains singular measures but without atoms. For item
(v), notice first that $\bar{\mu}$ is also a minimizer of $\Psi_{\nu}(\tau
,\mu_{0};\cdot).$ It follows that $\mu_{\nu,\tau}^{k}=$ $\bar{\mu}$ in
(\ref{scheme-1}), for all $k,$ and so $S[\bar{\mu}](t)=\mu_{\nu,\tau}%
(t)=\bar{\mu},$ for all $t>0$ and $\tau>0.$ Since $\mathcal{F}_{\nu}%
(\mu(t))\rightarrow\mathcal{F}_{\nu}(\bar{\mu})$ (by (\ref{func-conv})), as
$t\rightarrow\infty,$ the convergence $\mu(t)\rightarrow\bar{\mu}$ follows
from the compactness of $\mathcal{P}(\mathbb{S}^{1})$ and the lower
semicontinuity of $\mathcal{F}_{\nu}$. \fin

\begin{rem}
By item (iii), we can not exclude the possibility that in the case $\nu=0$ the
gradient flow of (\ref{eq-1}) is singular with respect to the Lebesgue
measure, although atoms in the solution are not allowed. In the viscous case,
the entropy part of the functional prevents the existence of singular measures
in the flow at $t>0$.
\end{rem}

\begin{teor}
\label{lim-inviscid}Let $\nu\geq0$ and, for $0<$ $\nu<\epsilon_{0},$ let us
denote by $\mu_{t}^{\nu}$ and $\mu_{t}$ the gradient flows associated to the
energy functionals $\mathcal{F}_{\nu}$ and $\mathcal{F}_{0}$, respectively,
with the same initial data $\mu_{0}\in\mathcal{P}(\mathbb{S}^{1})$. Then,
$\mu_{t}^{\nu}\rightarrow\mu_{t}$ in $\mathcal{P}(\mathbb{S}^{1})$, locally
uniformly in $[0,\infty)$, as $\nu\rightarrow0^{+}.$
\end{teor}

\emph{Proof.} By simplicity we can assume that $\mu_{0}\in D(\mathcal{F}%
_{\epsilon_{0}}).$ The general case follows by using the same argument
together with a discrete version of item (ii) of Theorem \ref{fluxograd} .

\textit{Step 1.} Let $\mu^{\nu}\rightarrow\mu$ as $\nu\rightarrow0^{+}$ and
let $\mu_{\tau}^{\nu}$ be the minimizer of $\Psi_{\nu}(\tau,\mu^{\nu};\cdot)$.
Then, for each fixed $\tau>0,$ $\mu_{\tau}^{\nu}\rightarrow\mu_{\tau}$ as
$\nu\rightarrow0^{+}$, where $\mu_{\tau}$ is the minimizer of $\Psi_{0}%
(\tau,\mu;\cdot)$.

In fact, by the compactness of $\mathcal{P}(\mathbb{S}^{1}),$ we can extract a
convergent subsequence $(\mu_{\tau}^{\nu_{l}})_{l=1}^{\infty}$ such that
$\nu_{l}\rightarrow0^{+}$ and $\rho=\lim_{l\rightarrow\infty}\mu_{\tau}%
^{\nu_{l}}.$ Given $\omega\in D(\mathcal{F}_{0})$, we have
\begin{equation}
\Psi_{\nu_{l}}(\tau,\mu^{\nu_{l}};\mu_{\tau}^{\nu_{l}})\leq\Psi_{\nu_{l}}%
(\tau,\mu^{\nu_{l}};\omega). \label{eq:liminviscid}%
\end{equation}
It is straightforward to check that $\lim_{l\rightarrow\infty}\Psi_{\nu_{l}%
}(\tau,\mu^{\nu_{l}};\omega)=\Psi_{0}(\tau,\mu;\omega)$. Also, by lower
semicontinuity, it follows that
\begin{align*}
\liminf_{l\rightarrow\infty}\mathcal{F}_{\nu_{l}}[\mu_{\tau}^{\nu_{l}}]  &
\geq\mathcal{F}_{0}[\rho]\\
\liminf_{l\rightarrow\infty}\mathbf{d}_{per}^{2}(\mu^{\nu_{l}},\mu_{\tau}%
^{\nu_{l}})  &  \geq\mathbf{d}_{per}^{2}(\mu,\rho).
\end{align*}
Then, we can conclude from (\ref{eq:liminviscid}) that $\rho$ is a minimizer
of $\Psi_{0}(\tau,\mu;\cdot)$ and so $\rho=\mu_{\tau}$ (by uniqueness).

\textit{Step 2.} Given $T>0$ finite and $\tau>0$, let $\mu_{\nu,\tau}(t)$ and
$\mu_{\tau}(t)$ be the discrete solution defined in (\ref{discrete-sol-1}). We
have that
\begin{equation}
\lim_{\nu\rightarrow0^{+}}\mathbf{d}_{per}(\mu_{\nu,\tau}(t),\mu_{\tau}(t))=0,
\label{eq:lim-inviscid1}%
\end{equation}
uniformly in $[0,T],$ for each $\tau>0$.

Recalling the definition of $\mu_{\tau}^{\nu,k}$ in (\ref{scheme-1}), notice
that \textit{Step 1} and an induction argument yield
\begin{equation}
\lim_{\nu\rightarrow0^{+}}\mathbf{d}_{per}(\mu_{\tau}^{\nu,k},\mu_{\tau}%
^{k})=0,\text{ for all }k\in\{0\}\cup\mathbb{N}. \label{aux-conv-1}%
\end{equation}
Since $\mu_{\nu,\tau}(t)$ and $\mu_{\tau}(t)$ are step functions with
$\mathcal{P}(\mathbb{S}^{1})$-values based on $\mu_{\tau}^{\nu,k}$ and
$\mu_{\tau}^{k}$ (see (\ref{discrete-sol-1})), respectively, the uniform
convergence (\ref{eq:lim-inviscid1}) in $[0,T]$ follows from (\ref{aux-conv-1}).

\textit{Step 3.} Finally, we can use Theorem \ref{fluxograd} (vi) in order to
get
\begin{align}
\mathbf{d}_{per}(\mu_{\nu}(t),\mu(t))  &  \leq\tau^{1/2}\left(  \mathcal{F}%
_{\nu}(\mu_{0})+2\pi\nu e^{-1}\right)  ^{1/2}+\mathbf{d}_{per}(\mu_{\nu,\tau
}(t),\mu_{\tau}(t))\label{aux-inviscid-1}\\
&  +\tau^{1/2}\mathcal{F}_{0}(\mu_{0})^{1/2},\text{ for all }t\in
\lbrack0,T].\nonumber
\end{align}
Since $\tau>0$ is arbitrary, we obtain the desired convergence by making
$\nu\rightarrow0^{+}$ in (\ref{aux-inviscid-1}) and using
(\ref{eq:lim-inviscid1}). \fin

\subsection{Solutions in distributional sense}

\hspace{0.4cm} We already know that the semigroup $\mu(t)=S[\mu_{0}](t)$ is an
absolutely continuous curve in $\mathcal{P}(\mathbb{S}^{1})$. The
characterization given in \cite{Ambrosio} for such curves provides that
$\mu(t)$ satisfies the continuity equation in the sense of distributions. More
precisely, for the periodic setting, the analogous one is given in
\cite[Theorem 3.2]{Gangbo}. In this case, the corresponding vector field is
the minimal selection in the sub-differential $\partial\mathcal{F}_{\nu}%
(\mu(t))$ and its $L^{2}$-norm with respect to $\mu(t)$ is a $L_{loc}^{1}$
function in $(0,\infty)$.

We say that a curve $\mu_{t}$ is a solution of the equation (\ref{eq-1}), with
initial data $\mu_{0}$, if, for each $\phi\in C^{\infty}(\mathbb{S}^{1}),$ the
equation
\begin{align}
\frac{d}{dt}\int_{-\pi}^{\pi}\phi(x)d\mu_{t}  &  =\nu\int_{-\pi}^{\pi}%
\phi^{\prime\prime}(x)d\mu_{t}\label{sol.dist}\\
&  -\int\!\int_{[-\pi,\pi)^{2}}\cot\left(  \frac{x-y}{2}\right)  \left(
\phi^{\prime}(x)-\phi^{\prime}(y)\right)  d(\mu_{t}\times\mu_{t})\nonumber
\end{align}
is verified in the distributional sense in $(0,\infty)$ and $\mu
_{t}\rightharpoonup\mu_{0}$ weakly-$\ast$ as measures, as $t\rightarrow0^{+}$.

For $\mu\in\mathcal{P}(\mathbb{S}^{1})$, we define $L_{\mu}:C^{\infty
}(\mathbb{S}^{1})\rightarrow\mathbb{R}$ by
\begin{equation}
L_{\mu}(\phi)=\int\!\int_{[-\pi,\pi)^{2}}\cot\left(  \frac{x-y}{2}\right)
(\phi(x)-\phi(y))d(\mu\times\mu). \label{energyfunctional-3}%
\end{equation}

The functional (\ref{energyfunctional-3}) is continuous in $C^{1}%
(\mathbb{S}^{1})$. In fact, since $\sin(0)/0=1$ (by definition), we can choose
$\delta>0$ such that $|t|<\delta$ implies that $\sin(t)/t\geq1/2$. Thus
\begin{align*}
|L_{\mu}(\phi)|  &  \leq\frac{1}{2\pi}\int\!\int_{[-\pi,\pi)^{2}\cap
\lbrack|x-y|<2\delta]}\left\vert \cot\left(  \frac{x-y}{2}\right)
(\phi(x)-\phi(y))\right\vert d(\mu\times\mu)\\
&  \ +\frac{1}{2\pi}\int\!\int_{[-\pi,\pi)^{2}\cap\lbrack|x-y|\geq2\delta
]}\left\vert \cot\left(  \frac{x-y}{2}\right)  (\phi(x)-\phi(y))\right\vert
d(\mu\times\mu)\\
&  \leq\frac{1}{\pi}[1+2\pi\cot(\delta)]\Vert\phi^{\prime}\Vert_{L^{\infty}}.
\end{align*}
Indeed, we can refine the above estimate by considering $\sin(t)/t\geq
\sin(\delta)/\delta$ for $|t|<\delta$. This shows the continuity of $L_{\mu
}:C^{\infty}(\mathbb{S}^{1})\rightarrow\mathbb{R}$. Now, arguing as in
\cite{Car-Fer-Pre} we have that there exists a Radon measure $\xi_{\mu}%
\in\mathcal{M}(\mathbb{S}^{1})$ such that the following representation for
(\ref{energyfunctional-3}) is valid
\begin{equation}
L_{\mu}(\phi)=\int_{[-\pi,\pi)}{\phi^{\prime}d\xi}_{\mu}{,}\text{ for all
}{\phi\in}C^{1}(\mathbb{S}^{1}). \label{funcw1}%
\end{equation}

\begin{lem}
\label{lemasubdif}For $\nu\geq0,$ let $\mathcal{F}_{\nu}$ be the functional
defined in (\ref{energyfunctional})-(\ref{energyfunctional1}) and $\mu\in
D(\mathcal{F}_{\nu})$. If $\mu\in D\left(  |\partial\mathcal{F}_{\nu}|\right)
$ then there exists $\theta\in L^{2}(\mathbb{S}^{1},d\mu)$ such that
\begin{equation}
\theta\mu=\frac{d}{dx}(\nu\mu+\xi_{\mu}), \label{subdif}%
\end{equation}
where $\xi_{\mu}$ is as in (\ref{funcw1}). The vector $\theta$ is the minimal
selection in the sub-differential $\partial\mathcal{F}_{\nu}(\mu)$. Moreover,
for $\nu>0$ we have $\nu\mu+\xi_{\mu}\in W^{1,1}(\mathbb{S}^{1},dx)$ where
$dx$ is the Lebesgue measure on $[-\pi,\pi)\equiv\mathbb{S}^{1}.$
\end{lem}

\emph{Proof.} We start by calculating the directional derivatives of
$\mathcal{F}_{0}$.

Let $\phi\in C^{\infty}(\mathbb{S}^{1})$, and define $r_{\epsilon
}(x)=x+\epsilon\phi(x)$. Note that $r_{\epsilon}(\pi)=2\pi+r_{\epsilon}(-\pi
)$, and so we can consider $r_{\epsilon}$ as a function from $\mathbb{S}^{1}$
to itself. Given $\mu\in D(\mathcal{F}_{\nu}),$ we can consider $r_{\epsilon
\#}\mu$ and, for $\epsilon<\text{Lip}_{\phi}^{-1},$ $r_{\epsilon
}(x)-r_{\epsilon}(y)>0$ provided that $x-y>0$. Analogously if $x-y<0$ then
$r_{\epsilon}(x)-r_{\epsilon}(y)<0$.

Therefore, since $-\frac{1}{\pi}\log\left\vert \sin\left(  \frac{x}{2}\right)
\right\vert $ is convex in each one of the segments $(0,2\pi]$ and
$[-2\pi,0),$ we can use monotonicity properties of convex functions in order
to obtain
\begin{align}
&  \frac{1}{\epsilon}\left[  \mathcal{F}_{0}(r_{\epsilon\#}\mu)-\mathcal{F}%
_{0}(\mu)\right] \nonumber\\
&  =-\frac{1}{\epsilon\pi}\int\!\int_{[-\pi.\pi)^{2}}\left(  \log\left\vert
\sin\left(  \frac{r_{\epsilon}(x)-r_{\epsilon}(y)}{2}\right)  \right\vert
-\log\left\vert \sin\left(  \frac{x-y}{2}\right)  \right\vert \right)
d\left(  \mu\times\mu\right)  . \label{aux-dist-1}%
\end{align}
Splitting the integral (\ref{aux-dist-1}) into the subsets $\{x>y\}$ and
$\{x<y\},$ we observe that each resulting integrand is nondecreasing in
$\epsilon$ (by convexity). By applying the monotone convergence theorem, we
obtain
\begin{align*}
\lim_{\epsilon\rightarrow0^{+}}\frac{1}{\epsilon}\left[  \mathcal{F}%
_{0}(r_{\epsilon\#}\mu)-\mathcal{F}_{0}(\mu)\right]   &  =-\frac{1}{2\pi}%
\int\!\int_{[-\pi,\pi)^{2}}\cot\left(  \frac{x-y}{2}\right)  \left(
\phi(x)-\phi(y)\right)  d\left(  \mu\times\mu\right) \\
&  =-L_{\mu}(\phi).
\end{align*}

Next we deal with the derivative of $\mathcal{U}$. Since $r_{\epsilon\#}\mu$
is supported in $[a,a+2\pi)$, for some $a\in\mathbb{R}$, we can use Lemma
\ref{entropyinvariance} to obtain
\[
\mathcal{U}[r_{\epsilon\#}\mu]=\int_{\mathbb{R}}f_{\epsilon}(x)\log
(f_{\epsilon}(x))\ dx,
\]
where $f_{\epsilon}$ is the density of $r_{\epsilon\#}\mu$. Now, the
computation of the directional derivative of $\mathcal{U}$ follows similarly
to that of the entropy functional in the whole space $\mathbb{R}$ (see
\cite{Villani}). It follows that
\begin{equation}
\lim_{\epsilon\rightarrow0}\frac{1}{\epsilon}[\mathcal{U}(r_{\epsilon\#}%
\mu)-\mathcal{U}(\mu)]=-\nu\int_{-\pi}^{\pi}{\phi^{\prime}d\mu}. \label{funcu}%
\end{equation}

We turn to the complete functional $\mathcal{F}_{\nu}$. First recall that
\[
\frac{\mathbf{d}_{per}(r_{\epsilon\#}\mu,\mu)}{\epsilon}\leq\Vert\phi
\Vert_{L^{2}(\mathbb{S}^{1},d\mu)}.
\]
Since $\mu\in D(|\partial\mathcal{F}_{\nu}|)$, we have that
\begin{align*}
\lim_{\epsilon\rightarrow0}\frac{1}{\epsilon}\left[  \mathcal{F}_{\nu
}[r_{\epsilon\#}\mu]-\mathcal{F}_{\nu}[\mu]\right]   &  =\nu\int_{-\pi}^{\pi
}{\phi^{\prime}d\mu}-L_{\mu}(\phi)\\
&  \geq-|\partial\mathcal{F}_{\nu}|(\mu)\Vert\phi\Vert_{L^{2}(\mathbb{S}%
^{1},d\mu)}.
\end{align*}
By changing $\phi$ by $-\phi$, it follows that
\[
\left\vert \nu\int_{-\pi}^{\pi}{\phi^{\prime}d\mu}-L_{\mu}(\phi)\right\vert
\leq|\partial\mathcal{F}_{\nu}|(\mu)\Vert\phi\Vert_{L^{2}(\mathbb{S}^{1}%
,d\mu)}.
\]
Therefore, by Riesz representation theorem, there exists $\theta\in
L^{2}(\mathbb{S}^{1},d\mu)$ such that
\begin{equation}
-\nu\int_{-\pi}^{\pi}{\phi^{\prime}d\mu}-\int_{[-\pi,\pi)}{\phi^{\prime}d\xi
}_{\mu}=\int_{-\pi}^{\pi}{\theta\phi d\mu,}\text{ for all }\phi\in C^{\infty
}(\mathbb{S}^{1}). \label{subdif1}%
\end{equation}
From (\ref{subdif1}), we have $\theta\mu=\frac{d}{dx}(\nu\mu+\xi_{\mu})$ in
the distributional sense and then $\theta$ is the minimal selection of
$\partial\mathcal{F}_{\nu}(\mu)$. Moreover, for $\nu>0,$ $\mu$ is absolutely
continuous with respect to the Lebesgue measure. Thus, by definition of
$\mathcal{U}$ and Holder inequality, it follows that $\theta\mu\in
L^{1}(\mathbb{S}^{1},dx)$, i.e. $\frac{d}{dx}(\nu\mu+\xi_{\mu})\in
L^{1}(\mathbb{S}^{1},dx)$, and thereby $\nu\mu+\xi_{\mu}\in W^{1,1}%
(\mathbb{S}^{1},dx)$. In particular, this implies that $\xi_{\mu}$ is
absolutely continuous with respect to the Lebesgue measure. \fin\bigskip

Now we are in position to show that our gradient flows are solutions for
(\ref{eq-1}) and (\ref{eq1-viscous}) in distributional sense.

\begin{teor}
\label{soldistt}Let $\nu\geq0$ and let $\mathcal{F}_{\nu}$ be the functional
defined in (\ref{energyfunctional})-(\ref{energyfunctional1}). Let $\mu_{t}$
be the gradient flow for $\mathcal{F}_{\nu}$ with initial data $\mu_{0}%
\in\mathcal{P}(\mathbb{S}^{1})$ and let $\xi_{\mu,t}$ be the Radon measure
associate to functional $L_{\mu_{t}}$ according to (\ref{funcw1}). Then,
$\mu_{t}:=\mu(x,t)$ is a weak solution of (\ref{eq1-viscous}) in the sense of
(\ref{sol.dist}). Moreover, for $\nu>0$ we have that
\begin{align}
\mu(x,t)  &  \in L_{loc}^{1}((0,\infty);W^{1,1}(\mathbb{S}^{1},dx))
\label{P1}\\
\nu\mu_{t}+\xi_{\mu,t}  &  \in W^{1,1}(\mathbb{S}^{1},dx),\text{ for all
}t>0\text{,} \label{P2}%
\end{align}
where $dx$ is the Lebesgue measure on $[-\pi,\pi)\equiv\mathbb{S}^{1}.$
\end{teor}

\emph{Proof.} We start by showing that $\mu_{t}$ is a distributional solution
of (\ref{eq1-viscous}). In fact, by definition of gradient flow, we know that
the continuity equation
\[
\partial_{t}\mu_{t}+\partial_{x}(\theta_{t}\mu_{t})=0
\]
is satisfied in the sense of distributions, where $\theta_{t}\mu_{t}=\frac
{d}{dx}(\nu\mu_{t}+\xi_{\mu,t}).$ For $\phi\in C^{\infty}(\mathbb{S}^{1})$ and
$\eta(t)\in C_{c}^{\infty}((0,\infty))$, we have that
\[
\int_{0}^{\infty}\int_{[-\pi,\pi)}[\eta^{\prime}(t)\phi(x)+\eta(t)\phi
^{\prime}(x)\theta_{t}(x)]\mu_{t}(x)dt=0
\]
and then%
\begin{align*}
&  \int_{0}^{\infty}\int_{[-\pi,\pi)}\eta^{\prime}(t)\phi(x)\mu_{t}(x){dt}\\
&  =\int_{0}^{\infty}\eta(t)\int_{[-\pi,\pi)}\phi^{\prime\prime}(x)d(\nu
\mu_{t}+\xi_{\mu,t})dt\\
&  =\nu\int_{0}^{\infty}\eta(t)\int_{[-\pi,\pi)}\phi^{\prime\prime}(x)d\mu
_{t}dt+\int_{0}^{\infty}\eta(t)L_{\mu_{t}}(\phi^{\prime})dt\\
&  =\nu\int_{0}^{\infty}\eta(t)\int_{[-\pi,\pi)}\phi^{\prime\prime}(x)d\mu
_{t}dt-\int_{0}^{\infty}\eta(t)\int\!\int_{[-\pi,\pi)^{2}}\cot(\frac{x-y}%
{2})(\phi^{\prime}(x)-\phi^{\prime}(y))d\left(  {{\mu_{t}\times\mu}}%
_{t}\right)  dt,
\end{align*}
which gives (\ref{sol.dist}).

Next, for $\nu>0,$ notice that
\[
\lVert\partial_{x}\mu(\cdot,t)\rVert_{L^{1}(\mathbb{S}^{1},dx)}=\lVert
\frac{\partial_{x}\mu(x,t)}{\mu(x,t)}\rVert_{L^{1}(\mathbb{S}^{1},d\mu_{t}%
)}\leq\lVert\frac{\partial_{x}\mu(x,t)}{\mu(x,t)}\rVert_{L^{2}(\mathbb{S}%
^{1},d\mu_{t})}.
\]
The properties (\ref{P1})-(\ref{P2}) follows from Lemma \ref{lemasubdif} and
the fact that $\frac{\partial_{x}\mu(x,t)}{\mu(x,t)}$ is the minimal selection
for the functional $\mathcal{U}$ in $\mu_{t}$, and then $\lVert\frac
{\partial_{x}\mu(x,t)}{\mu(x,t)}\rVert_{L^{2}(\mathbb{S}^{1},d\mu_{t})}\in
L_{loc}^{1}((0,\infty))$. \fin

\end{document}